\documentclass[11pt]{article}

\usepackage{comment}
\usepackage{amsmath,amssymb,amsthm,comment}
\usepackage[colorlinks,citecolor=red,urlcolor=blue,linkcolor=blue]{hyperref}
\usepackage{mathrsfs,dsfont}
\usepackage{graphicx}
\usepackage{enumitem}
\usepackage{nicefrac}
\usepackage[font={normalsize},width=1.12\linewidth]{caption}

\usepackage{tikz,graphicx}
\usetikzlibrary{calc,arrows.meta}

\newtheorem{nummer}{ }
\newtheorem{thm}[nummer]{\bf Theorem}

\newtheorem{cor}[nummer]{\bf Corollary}

\newtheorem{exa}[]{\bf Example}

\newcommand{\ie} {i.e.}

\newcommand{\equimod}[2]{{#1}\,(\operatorname{mod}\,{#2})}

\renewcommand{\phi}{\varphi}
\renewcommand{\theta}{\vartheta}

\newcommand{\Q}{\mathds{Q}}

\newcommand{\Z}{\mathds{Z}}

\newcommand{\C}{\mathds{C}}
\newcommand{\F}{\mathds{F}}

\newcommand {\mult}{\mathbin{*}}

\makeatletter
\def\opargproof[#1]{\par\noindent {\bf #1 }}
 \makeatother

\begin{document}
\begin{center}
{\LARGE\bf Halving formulae for points on elliptic curves}

\bigskip
{\small Lorenz Halbeisen}\\[1.2ex] 
{\scriptsize Department of Mathematics, ETH Zentrum,
R\"amistrasse\;101, 8092 Z\"urich, Switzerland\\ lorenz.halbeisen@math.ethz.ch}\\[1.8ex]

{\small Norbert Hungerb\"uhler}\\[1.2ex] 
{\scriptsize Department of Mathematics, ETH Zentrum,
R\"amistrasse\;101, 8092 Z\"urich, Switzerland\\ norbert.hungerbuehler@math.ethz.ch}
\end{center}

\hspace{5ex}{\small{\it key-words\/}: Elliptic curves, halving formulae, finite fields}

\hspace{5ex}{\small{\it 2020 Mathematics Subject 
Classification\/}: {\bf 11G05}\,\ 94A60\ 11R16}

\begin{abstract}\noindent
Let $P$ be an arbitrary point on an elliptic curve 
over the complex numbers of the form 
\hbox{$y^2=x^3+a_4\,x+a_6$} or of the form
\hbox{$y^2=x^3+a_2\,x^2+a_4\,x$}.
We provide explicit formulae to compute the points $P/2$,
{\ie}, the points $Q$ such that $2\mult Q=P$. 
\end{abstract}

\section{Introduction}

Even though it is well-known how to double point on elliptic 
curves, very little is known about halving points. One of the
few results, which we obtain as Corollary\;\ref{cor:e012}, 
is that a rational point $P=(x_0,y_0)$ on an elliptic curve 
$E:\;y^2=(x-e_0)(x-e_1)(x-e_2)$ (with $e_i\in\Q$) is in
$2E(\Q)$, if and only if $x_0-e_0$, $x_0-e_1$, $x_0-e_2$ are all
squares of rational numbers
(see \cite[Ch.\,I,\,\S9,\,Prp.\,20]{Koblitz}). 
For related results we refer the reader to 
\cite{Bekker,Sadek,Yelton,Zarkhin}.

The natural approach to finding formulas for halving points
is to invert the doubling formula. However, this leads in
general to quartic equations which are hard to solve. Another
approach is to extend first the basic field in order to make
the equations bi-quadratic, which is what we do below.

In the next section, we provide explicit formulae for halving
points on elliptic curves over~$\C$ of the form
\hbox{$y^2=x^3+a_4\,x+a_6$} and of the form
\hbox{$y^2=x^3+a_2\,x^2+a_4\,x$}, respectively. 
Notice that since every elliptic
curve can be transformed into this form, similar formulae 
for halving points exist for arbitrary elliptic curves. Moreover,
the formulae work also if we replace $\C$ by some finite field.
As an example we show how to halve points on the curves 
$y^2=x^3-n^2+x$ for positive integers~$n$ and on curves
of the form $y^2=x^3+a_2x+a_6$ over prime fields.

\section{Formulae for $\boldsymbol{P/2}$}

\begin{thm}\label{thm:main}
Let\/ $E_{a_4,a_6}$ and\/ $E_{a_2,a_4}$ 
be non-singular elliptic curves over\/ $\C$ defined by
$$E_{a_4,a_6}:\ y^2 = x^3 + a_4\,x + a_6\qquad\text{and}\qquad
E_{a_2,a_4}:\ y^2 = x^3 + a_2\,x^2 + a_4\,x$$
where\/ $a_2,\,a_4,\,a_6\in\C$, 
and let\/ $P=(x_0,y_0)$ be a point on~$E_{a_4,a_6}$ or 
on~$E_{a_2,a_4}$.

If\/ $P\in E_{a_4,a_6}$, let
$$
r = -9\,a_6 + \sqrt{3}\,\sqrt{4\,a_4^3 + 27\,a_6^2},\qquad
d =-a_4\,\Bigl(\frac{2}{3r}\Bigr)^{\frac 13}+
\Bigl(\frac{r}{18}\Bigr)^{\frac 13},\qquad
k = -(a_4 + 3d^2),
$$
and
$$
e_0 = d,\qquad
e_1 = \frac{-d+\sqrt{9d^2+4k}}{2},\qquad
e_2 = \frac{-d-\sqrt{9d^2+4k}}{2}.\qquad
$$
If\/ $P\in E_{a_2,a_4}$, let
$$e_0=0,\qquad
e_1=\frac{-a_2 + \sqrt{a_2^2 - 4 a_4}}2,\qquad
e_2 =\frac{-a_2 - \sqrt{a_2^2 - 4 a_4}}2.
$$
Finally, in both cases let 
$$\gamma =\sqrt{x_0-e_0},\qquad
\alpha =\sqrt{x_0-e_1},\qquad
\beta  =\sqrt{x_0-e_2}.$$
Then, both elliptic curves\/ $E_{a_4,a_6}$ and\/ $E_{a_2,a_4}$ are of the
form
$$y^2=(x_0-e_0)(x_0-e_1)(x_0-e_2)$$
and the\/ $x$-coordinates of the four points\/ 
$Q_1,Q_2,Q_3,Q_4$ with\/ $2*Q_i=P$ are
\begin{eqnarray*}
x_{11} & = & x_0 + \alpha\beta + \gamma(\alpha+\beta),\\[.8ex]
x_{12} & = & x_0 + \alpha\beta - \gamma(\alpha+\beta),\\[.8ex]
x_{21} & = & x_0 - \alpha\beta + \gamma(\alpha-\beta),\\[.8ex]
x_{22} & = & x_0 - \alpha\beta - \gamma(\alpha-\beta).
\end{eqnarray*}
\end{thm}

\begin{proof} With respect to the curve $E_{a_2,a_4}$
we obviously have $-e_0-e_1-e_2=a_2$, $e_0e_1+e_1e_2+e_2e_0=a_4$,
and $-e_0e_1e_2=0$.

With respect to the curve $E_{a_4,a_6}$
we have to show that $-e_0-e_1-e_2=0$, which is obvious,
that $e_0e_1+e_1e_2+e_2e_0=a_4$, and that $-e_0e_1e_2=a_6$.

We have $$e_0\,e_1+e_1\,e_2+e_2\,e_0\;=\;-d^2+
\frac{d^2-9d^2-4k}{4}=-3\,d^2-k$$
and $$-e_0\,e_1\,e_2\;=\;-d\,\frac{d^2-9d^2-4k}{4}
=d(2\,d^2+k).$$
So, we have to show that 
$$a_4 = -(3\,d^2+k)\qquad\text{and}\qquad
a_6 = d(2\,d^2+k).$$ 
By definition of $k$ we have
$$-(3d^2+k)=-(3d^2-a_4-3d^2)=a_4$$
and 
\begin{multline*}
d(2d^2+k)=2d^3-d\,a_4-3d^3=-d^3-d\,a_4=\\[1ex]
\biggl(a_4^3\,\Bigl(\frac{2}{3r}\Bigr)
\underset{=da_4}{\underbrace{
- a_4^2\,\Bigl(\frac{2}{3r}\Bigr)^{\frac 13}+
a_4\,\Bigl(\frac{r}{18}\Bigr)^{\frac 13}}}
-\frac{r}{18}\biggr)-d\,a_4 = a_4^3\,\frac{2}{3r}-
\frac{r}{18}=\qquad\qquad\\[1ex]
\frac{12\,a_4^3-r^2}{18\,r}=
\frac{12\,a_4^3-81\,a_6^2-12\,a_4^3 -81\,a_6^2+
18a_6\,\sqrt{4\,a_4^3 + 27\,a_6^2}}
{-162\,a_6+18\,\sqrt{4\,a_4^3 + 27\,a_6^2}}=\\[1ex]
\frac{a_6\bigl(-162\,a_6+18\,\sqrt{4\,a_4^3 + 27\,a_6^2}\,\bigr)}
{-162\,a_6+18\,\sqrt{4\,a_4^3 + 27\,a_6^2}}
=a_6.
\end{multline*}
By shifting the curve, we may assume that $e_0=0$. In particular,
we may assume that the elliptic curve is $E_{a_2,a_4}$.
To show that $x_{11},x_{12},x_{21},x_{22}$ are the 
the $x$-coordinates of points $Q\in E_{a_2,a_4}$ such that $2\mult Q=P$,
it is enough to show that the $x$-coordinate 
of the point $Q_{ij}:=(x_{ij},y)$, 
where $i,j\in\{1,2\}$ and $y=\sqrt{x_{ij}^3+a_2\,x_{11}^2+a_4\,x_{11}}$,
is equal to~$x_0$. Now, the $x$-coordinate $x_{2ij}$
of the point $2\mult Q_{ij}$
is given by the formula
$$x_{2ij}=\frac{x_{ij}^4-2a_4\,x_{ij}^2+a_4^2}
{4(x_{ij}^3+a_2\,x_{11}^2+a_4\,x_{11})}=\frac{(x_{ij}^2-a_4)^2}
{4x_{ij}(x_{ij}^2+a_2\,x_{ij}+a_4)}.$$
Now, if $e_0=0$, then $x_0=\gamma^2$ and we obtain
\begin{eqnarray*}
x_{11} & = & (\alpha+\gamma)(\beta+\gamma),\\[.8ex]
x_{12} & = & (\alpha-\gamma)(\beta-\gamma),\\[.8ex]
x_{21} & = & (\alpha+\gamma)(-\beta+\gamma),\\[.8ex]
x_{22} & = & (\alpha-\gamma)(-\beta-\gamma).
\end{eqnarray*}
Furthermore, we have $a_2=\alpha^2 + \beta^2 - 2\gamma^2$ 
and $a_6=(\alpha^2 - \gamma^2) (\beta^2 - \gamma^2)$, and if we 
write $x_{ij},a_2,a_4$ in terms of $\gamma,\alpha,\beta$, it is
not hard to verify that 
$$(x_{ij}^2-a_4)^2=4x_{ij}\gamma^2(x_{ij}^2+a_2\,x_{ij}+a_4),$$
which shows that $x_{2ij}=x_0$.
\end{proof}

\noindent{\it Remarks.} (a)~As a matter of fact 
we would like to mention that 
the lines through the pairs of 
points with $x$-coordinate $x_{11},x_{12}$ and $x_{21},x_{22}$, 
respectively, meet in a point $S=(x_S,y_S)$ on the curve~$E$ with
$$x_S= 
d+\frac{k}{d-x_0}
\hspace{1.2ex}\text{for $E=E_{a_4,a_6}$}\quad\text{and}\quad
x_S=\frac{a_4}{x_0}\hspace{1.2ex}\text{for $E=E_{a_2,a_4}$ ({\ie},
$d=0,\,k=-a_4$).}$$

(b)~If $E$ is an elliptic curve of the form $E_{a_4,a_6}$ or
$E_{a_2,a_4}$ over the field $\F$ and $P,Q\in E(\F)-0$ are such that
$2\mult Q=P$, then $Q$ is unique if and only if there is no element
in the group $(E(\F),+)$ of order~$2$. 

(c)~If we set $t = x_0^2+dx_0-2d^2-k$ and
$$
w =\sqrt{t},\qquad
w_1= \sqrt{(x_0-d)(d+2w+2x_0)},\qquad
w_2=\sqrt{(x_0-d)(d-2w+2x_0)},$$
then, for $E_{a_4,a_6}$, we have 
$$x_{11} = x_0 + w + w_1,\quad
x_{12} = x_0 + w - w_1,\quad
x_{21} = x_0 - w + w_2,\quad
x_{22} = x_0 - w - w_2.
$$

\smallskip

As a consequence of Theorem\;\ref{thm:main} we obtain the following

\begin{cor}\label{cor:e012}
A rational point\/ $P=(x_0,y_0)$ on an elliptic curve
$$E:\;y^2=(x-e_0)(x-e_1)(x-e_2)\quad
\text{with\/ $e_0,e_1,e_2\in\Q$}$$ is 
equal to\/ $2\mult Q$ for some rational point~$Q$ on~$E$,
if and only if\/ $x_0-e_0$, $x_0-e_1$, $x_0-e_2$ are all
squares of rational numbers.
\end{cor}

\begin{proof}
Notice first that since the torsion group of $E(\Q)$ is isomorphic
to $\Z_2\times\Z_2$, for any rational point $P=(x_0,y_0)$ on~$E$,
either none or all four points $Q_i$ with $2\mult Q_i$ are rational.
Now, since we have $y_0=-\lambda x_0+\lambda x_{ij}-y_{ij}$, where
$$\lambda = \frac{3x_{ij}^2+2a_2 x_{ij}+a_4}{2y_{ij}},$$
and $y_0$ is rational, $y_{ij}$ is rational if and only if 
$x_{ij}$ is rational. In particular, either none or all four 
values $x_{11},x_{12},x_{21},x_{22}$ are rational. So, it is enough
to show that $x_{11}$ is rational if and only if 
$x_0-e_0$, $x_0-e_1$, $x_0-e_2$ are all
squares of rational numbers.
 
By shifting the curve $E$, we may assume that $e_0=0$, which gives us
$$x_{11}= (\alpha+\gamma)(\beta + \gamma).$$ 
\noindent $(\Leftarrow)$ 
If $x_0=\gamma^2$, $x_0-e_1=\alpha^2$, $x_0-e_2=\beta^2$ 
are squares of rational numbers $\gamma,\alpha,\beta$, respectively,
then $x_{11}$ is obviously a rational number.

\noindent $(\Rightarrow)$ 
If $x_{11},x_{12},x_{21},x_{22}$ are rational, then all
$$
(\alpha+\gamma)(\beta + \gamma),\quad
(\alpha-\gamma)(\beta - \gamma),\quad
(\alpha+\gamma)(-\beta + \gamma),\quad
(\alpha-\gamma)(-\beta - \gamma),$$
are rational, which happens just in the case when 
$\gamma,\alpha,\beta$ are all rational.
\end{proof}

\section{Two Examples}

\begin{exa}{\rm
For natural numbers $n$, let us consider the curve 
$$E_n:\ y^2=x^3-n^2x$$ over the field~$\Q$.
It is well known that $n$ is a congruent number if and only if 
there is a rational point $(x_0,y_0)$ on $E_n$ with $y_0\neq 0$. 
Further it is known that the torsion group of $E_n(\Q)$ is isomorphic
to $\Z_2\times\Z_2$\,---\,notice that the three points of order~$2$ are
$(0,0)$ and $(\pm n,0)$. 

Now, for $a_4=-n^2$ and $a_6=0$, we have $e_0=0$, $e_1=n$, and $e_2=-n$.
Assume that $(x_0,y_0)$ is a rational point on~$E$ with $y_0\neq 0$.
In particular, $(x_0,y_0)$ is a point of infinite order.
Then our halving formulae give us 
\begin{eqnarray*}
x_{11} & = & x_0 + 
\sqrt{x_0^2-n^2} + \sqrt{x_0}\,\bigl(\sqrt{x_0-n}+\sqrt{x_0+n}\bigr),\\[1ex]
x_{12} & = & x_0 + 
\sqrt{x_0^2-n^2} - \sqrt{x_0}\,\bigl(\sqrt{x_0-n}+\sqrt{x_0+n}\bigr),\\[1ex]
x_{21} & = & x_0 - 
\sqrt{x_0^2-n^2} + \sqrt{x_0}\,\bigl(\sqrt{x_0-n}-\sqrt{x_0+n}\bigr),\\[1ex]
x_{22} & = & x_0 -
\sqrt{x_0^2-n^2} - \sqrt{x_0}\,\bigl(\sqrt{x_0-n}-\sqrt{x_0+n}\bigr).
\end{eqnarray*}
For $n=6$, $x_0=-3$, and $y_0=9$, the point $P=(x_0,y_0)$ is a point of
$E_n$ of infinite order and since $\sqrt{x_0^2-n^2}\notin\Q$, 
there is no rational point $Q\in E_n$ such that $2\mult Q=P$.
On the other hand, $2\mult P=(25/4,-35/8)$, and for the four points
$$Q_1=(18,-72),\quad Q_1=(-2,-8),\quad Q_1=(-3,9),\quad Q_1=(12,36),$$
we have that $2\mult Q_i=2\mult P$.
}
\end{exa}

\begin{exa}{\rm
Let $p$ be an odd prime with
$p>\max\{|a_4|,|a_6|\}$, and let $E$ be an elliptic curve
of the form
$$E:\ y^2=x^3+a_4\,x+a_6$$ over the field~$\F_p$.
Furthermore, let $P=(x_0,y_0)$ be a point on~$E(\F_p)$.

Assume that the polynomial $X^3+a_4\,X+a_6\in\F_p[X]$ is 
irreducible over $\F_p$, which implies that no element of
the group $(E,+)$ has order~$2$.  
Then there exists a unique $Q\in E(\F_p)$
such that $Q=P/2$, {\ie}, $2\mult Q=P$. Now, consider the field 
$\F_{p^3}:=\F_p[X]/(X^3+a_4\,X+a_6)$. 
Then, for $d=X$ and $k=-(a_4+3\,X^2)$, the calculations to compute 
$x_{ij}$ (for $i,j\in\{1,2\}$) take place in the field~$\F_{p^3}$,
where exactly on of $x_{ij}$ belongs to $\F_{p^3}$, which is the
$x$-coordinate of $P/2$. We would like mention that
in the case when $p\equiv\equimod 34$, 
for any $a\in\F_{p^3}$ which is
a square we have $\sqrt{a}=\pm a^{\nicefrac{(p^3+1)}{4}}$ (for the 
computations of roots in the case when $p\equiv\equimod 34$ 
see~\cite{Mueller}).

To illustrate the algorithm, let 
$$p=17000000000000071,\quad
E:\ y^2=x^3+17\,x+71,\quad
P=(17071,\;4145148307074498).$$
Then $P\in E(\F_p)$, $p\equiv\equimod 34$, 
and the polynomial $f=X^3+17X+71$ is
irreducible over $\F_p$. In order to compute $P/2$, we
work in the field $\F_{p^3}=\F_p[X]/(f)$ and compute first
$\pm w$, $w_1$ and $w_2$\,---\,given in 
Remark\;(c)\,---\,with respect 
to $x_0=1700000000000071$, $d=X$, and $k=-(17+3X^2)$:
\begin{eqnarray*}
w\phantom{{}_1} & = & 
14551045313109763\,X^2 + \phantom{1}6554553633085449\,X + 
7788921359847751\\
-w\phantom{{}_1} & = &
\phantom{1}2448954686890308\,X^2 + 
10445446366914622\,X + 9211078640152320\\
w_1 & = &
\phantom{1}2448954686890308\,X^2 + 
10445446366914622\,X + 9211078640152320\\
w_2 & = &
11464660070096309\,X^2 + \phantom{1}3152403672687182\,X + 
16599480794424362.
\end{eqnarray*}
Thus, $w+w_1=w_1-(-w)\in\F_p$ and since $x_0\in\F_p$ we have
$$x_{11}=x_0+w+w_1=4631223433830370\in\F_p$$
which gives us the point
$$P/2=(4631223433830370,\;13664114850453464),$$
and since
$x_{12},x_{21},x_{22}\notin\F_p$, $P/2$ is the unique 
point on $E(\F_p)$ with the property that $2\mult(P/2)=P$.

As a matter of fact we would like to mention that we obtain the 
point $P/2$ also by first computing the order of $P$ in the 
finite group $(E,+)$: 
The order of $P$ is $m=16999999816127027$ (which is the same as the 
order of~$E$), and since the multiplicative inverse of~$2$ in $\Z_m$ 
is $\bar 2=8499999908063514$, we have $P/2=\bar 2\mult P$.
}
\end{exa}

\noindent{\it Remark.} By our halving formulae we can invert the
Diffi-Hellman key exchange algorithm with respect to some elliptic 
curves~$E$ over a finite field and a point $P$ on~$E$: 
If the order of the elliptic curve~$E$ is odd, then, since
the group $(E,+)$ is a finitely generated abelian group, 
for each point $P$ on $E$, $P/2$ is unique. Now,
after Alice and Bob have exchanged a common secret key $k$, 
Alice can encrypt a positive integer $T$ as follows. First, Alice
multiplies $T$ by some integer $10^m$ and searches a point 
$Q=(10^m\cdot T+i,y_i)$ on $E$ where $i<10^m$. Then, Alice goes
step by step through the binary representation of $k$, and in
each step, beginning with $Q$, she doubles the point and add the
point $b_n\mult P$, where $b_n$ is the $n$th entry in the binary
representation of~$k$. Alice sends the result $Q'$ to Bob, who can 
invert Alice's algorithm, starting with $Q'$,
by reading the binary representation of the key $k$ 
from the back, by replacing addition with subtraction of $P$, 
and by replacing doubling by halving points.

\bibliographystyle{plain}

\end{document}